\documentclass[a4paper, 11pt]{article}

\usepackage[english]{babel}
\usepackage{graphicx} \graphicspath{{./}{Figures/}}
\usepackage{amsmath}
\usepackage{amsfonts}
\usepackage{array}
\usepackage{subfig}
\usepackage{microtype}
\usepackage{epstopdf}
\usepackage{a4wide}

\usepackage{algorithm}
\usepackage{algorithmicx}
\usepackage{algpseudocode}

\usepackage{color}

\def\Prob{\mathbb{P}}
\def\E{\mathbb{E}}
\def\Var{\mathrm{Var}}

\newcommand\blfootnote[1]{%
  \begingroup
  \renewcommand\thefootnote{}\footnote{#1}%
  \addtocounter{footnote}{-1}%
  \endgroup
}

\title{The age of information in gossip networks}

\author{Jori Selen\footnotemark[2], Yoni Nazarathy\footnotemark[3] \footnotemark[4], Lachlan L.H. Andrew\footnotemark[4], and Hai L. Vu\footnotemark[4]}

\begin{document}
\maketitle

\renewcommand{\thefootnote}{\fnsymbol{footnote}}
\footnotetext[2]{Department of Mathematics and Computer Science, Eindhoven University of Technology, The Netherlands}
\footnotetext[3]{School of Mathematics and Physics, The University of Queensland, Australia}
\footnotetext[4]{Faculty of ICT, Swinburne University of Technology, Australia}
\blfootnote{E-mail adresses: {\tt j.selen@tue.nl}, {\tt y.nazarathy@uq.edu.au}, {\tt \{landrew,hvu\}@swin.edu.au}}
\renewcommand{\thefootnote}{\arabic{footnote}}

\begin{abstract}
We introduce models of gossip based communication networks in which each node is simultaneously a sensor, a relay and a user of information. We model the status of ages of information between nodes as a discrete time Markov chain. In this setting a gossip transmission policy is a decision made at each node regarding what type of information to relay at any given time (if any). When transmission policies are based on random decisions, we are able to analyze the age of information in certain illustrative structured examples either by means of an explicit analysis, an algorithm or asymptotic approximations. Our key contribution is presenting this class of models.
\end{abstract}

\section{Introduction} \label{sec:introduction}
We consider gossip networks in which the nodes wish to maintain an updated situation awareness view of the information sensed by all other nodes in the network. Using the {\em gossip} paradigm \cite{demers,karp}, this is done by having nodes transmit both their own sensed information and information that they have received from others. Thus nodes act as sensors, relays and receivers. Bandwidth is limited and communication channels are imperfect, thus the decision of what and when to transmit may often greatly affect performance. A natural application for gossip networks is intelligent transport systems (ITS) in which vehicles wirelessly share information relating to traffic congestion, road conditions and route alternatives, in order to improve safety and reduce congestion \cite{dimitrakopoulos,papadimitratos}. In this setting, gossiping is a suitable way to overcome the frequent changes in network topology.

The decision at each node of whether to transmit and what to transmit, are typically taken so as to minimize some measure of cost. Natural measures include the {\em ages of information} between the various node pairs, where the age of information at node $i$ of information sensed at node $j$ is defined as the difference between the current time and the time-stamp found on the most recent sensor measurement from $j$ received (perhaps through relays) at $i$.

Our aim is to introduce simple Markovian age of information models together with preliminary performance analysis results. Such models may influence network planning, protocol design and synthesis of efficient control methods. For the specific examples in this paper, it is easy to generate efficient deterministic transmission policies, but the analysis we carry out here is a first step toward studying more complex networks in which randomized policies are beneficial.

A fundamental question in the design of gossip networks is the following: {\em In order to help the greater good, how should a node balance relaying with transmitting its own information?} This paper sets the tone for treatment of this question by means of performance analysis and optimal policy design. For the specific case of ring networks, we give an answer based on asymptotics.

There has been much work focusing on either information aggregation
\cite{boyd,jelasity,kempe} or the age of information in gossip networks
\cite{bakhshi,banjeree,chaintreau,eugster}. The former dealt with the problem
of computing aggregates based on some functions, such as sum, average or
quantile of a set of data distributed over the nodes of a gossip network, and
studied the performance of protocols in terms of convergence and the optimization of neighbour selection (i.e.~strategy). The latter looked at the age of information via either analyzing the evolution of processes that gossip one message or content \cite{banjeree,eugster} or characterizing the distribution of latency (i.e.~age) over the network of many nodes \cite{bakhshi,chaintreau}.

In particular, both models in \cite{bakhshi,chaintreau} are based on a mean
field analysis with the networks size tending to infinity. The model in
\cite{chaintreau} yields a set of partial differential equations that uniquely
describe a system that allowed opportunistic content updates as in our work but
without interference or a lossy wireless channel. The model in \cite{bakhshi},
on the other hand, is based on a discrete-time Markov chain which could
possibly be extended to account for a lossy channel but without a content
update. Finally, \cite{fehnker} considers a lossy channel, and uses model
checking and Monte Carlo simulation to investigate the performance of a probabilistic broadcast gossip protocol.

Asymptotic results for a problem related to the age of information have been studied under the name of first-passage percolation~\cite{hammersley1965first}. Results in that field typically consider a single piece of information spreading on an infinite two dimensional lattice, and consider properties such as the shape of the region which has obtained information by a given time~\cite{smythe1978first}, or the variance of the time until the information reaches a given location~\cite{benjamini2003first}. Much
less work has considered irregular networks, although there has been some study of the Dirichlet triangulation of a two-dimensional Poisson process~\cite{vahidi1990first},  geometric graph networks $\mathbb{R}^d$~\cite{friedrich2011diameter}, Erdosh-Renyi networks~\cite{van2001first} and scale-free networks~\cite{bhamidi2010extreme}.

Our models and flavour of results are different in that we propose a simpler Markovian framework that can provide explicit formulae for the stationary distribution of the age of information in some specific cases. Using this framework the mean age at each node is also obtained for arbitrary tree networks, while the same is achieved via asymptotic analysis for ring networks. A further distinctive feature is that our models are suited to real-time data that is continuously updated. This differs from models where one big file is being transferred, or sensor network models where the key aim is to conserve energy, as in \cite{goseling}.

%
The rest of this paper is organized as follows.
Section~\ref{sec:modelDescription} introduces the age of information models.
These are specialized to linear, tree and ring networks in
Section~\ref{sec:structModels}, where we also present some basic results for the
mean and variance of the age of information and motivate the understanding of
rings. Section~\ref{sec:explicit} presents some non-trivial explicit and
algorithmic solutions for specific structured examples.
Section~\ref{sec:asymptoticsAndOptimPolicy} presents asymptotic approximations
for structured ring networks with a simple policy where we also answer the
question of the balance between relaying and transmitting one's own information.

\section{Age of information modeling} \label{sec:modelDescription}
We consider networks of a finite number of nodes, in which sensing, transmission and reception occurs at discrete (slotted) time instances. The age of information process, $\{A_{i,j}(n),~ n=0,1,2,\ldots\}$ is such that $A_{i,j}(n)$ is the age of the information that node $i$ has about node $j$ at time $n$. Thus for example if $A_{1,3}(n)=15$, we know that at time $n$, node $1$'s most updated view regarding the sensed information at node $3$ is from time $n-15$.

We denote the sequence of information transmissions indicators by $\{I_{i,j}(n)\}$, where $I_{i,j}(n) = 1$ if and only if at time $n$ node $i$ has broadcast its information regarding node $j$, otherwise $I_{i,j}(n)=0$. Note that $I_{i,i}(n)$ indicates if a node broadcasts its own sensed information.

We assume some sort of channel model in which the received packets at a given time $n$ at every node depend on the transmitted packets in the whole network at time $n$ and some other possible random effects that are independent for different $n$, yet follow the same probabilistic law. This may describe essentially any form of time-independent communication channel without memory.
At time $n$ the resulting receptions of packets are a random function of
$I_{i,j}(n)$ for all $i,j$ and are denoted by $R_{i,j}(n)$ where $R_{i,j}(n)=1$
if and only if $j$ received a packet sent by node $i$ (containing any form of
sensor information, original or relayed). Using $\wedge$ to denote the minimum,
the dynamics of the age process are
\begin{align}
\label{eqn:mainRecurs}
A_{i,j}(n+1) &= \left\{ \begin{array}{ll}
\left( A_{i,j}(n)~ \wedge~ \bigwedge_{\{k: R_{k,i}(n) I_{k,j}(n) =1\}} A_{k,j}(n)  \right) + 1, & i \neq j \, , \\
0, & i=j \, .
\end{array} \right. &
\end{align}
As \eqref{eqn:mainRecurs} illustrates, age increases by $1$  at each time slot, unless ``fresh information'' is received. Each node $i$ is only interested in the ``freshest'' information about $j$ and therefore compares the minimum age of information that was received (on all receptions $k$) with the current age of information stored in node $i$. The channel plays a role here in determining how $I(n)$ is mapped to $R(n)$: $I(n)$ determines all transmissions made on the network and this in turn (perhaps taking interference into account) determines all receptions.

Randomness enters \eqref{eqn:mainRecurs} through both the channel and possibly
through the transmission decisions $I(n)$ in case they are random. In this
paper we shall take $\{I(n)\}$ to be a (multi-dimensional) i.i.d. sequence. We
refer to this as having {\em Bernoulli policies}, i.e., the decision of what to transmit at any time instant is based on the time-invariant probability distribution of $I(n)$. In this case it is clear that \eqref{eqn:mainRecurs} together with some initial distribution, defines a discrete time Markov chain.

For a network of $N$ nodes where each node is assumed to have a sensor, the
state space of the Markov chain is $\mathbb{Z}_+^{N^2-N}$. Transitions on this
space are either of the form (a)~incrementing a coordinate by $1$ (no new
reception) or (b)~shifting a coordinate to equal the value of another
coordinate plus $1$ (new reception of fresh information). Showing that the
Markov chain has a single irreducible countably infinite class (nicely
represented as a subset of $\mathbb{Z}_+^{N^2-N}$), is non-periodic and is
positive recurrent, is straight-forward under quite general assumptions on the
channels and the transmission policy. We shall skip these details as they are
non-instructive. (Positive recurrence can be established by means of a linear Lyapunov function.)

Finding explicit performance measures, most importantly finding the stationary distribution, marginals of the stationary distribution or their mean, poses a much greater challenge. In the remainder of the paper we focus on introductory special structured examples on which the behaviour can by analyzed.

\section{Structured models}
\label{sec:structModels}
In order to get some insight into the behaviour of age of information models of the form \eqref{eqn:mainRecurs}, we look at some structured examples. To do so, we assume that the channel is represented by a directed graph, indicating which nodes can directly communicate. The graph determines the possible paths in which information may flow from sensor to user. The minimal attainable age of information, $A_{i,j}(n)$, is then the shortest path on the graph from $j$ to $i$. In case there is no such path, the $A_{i,j}(\cdot)$ component of the Markov chain is ignored.

\subsection*{Linear and tree networks}
\begin{figure}
\vskip-4ex
\subfloat[]{%
\begin{minipage}[c][0.5\width]{0.5\textwidth}
\centering%
\includegraphics{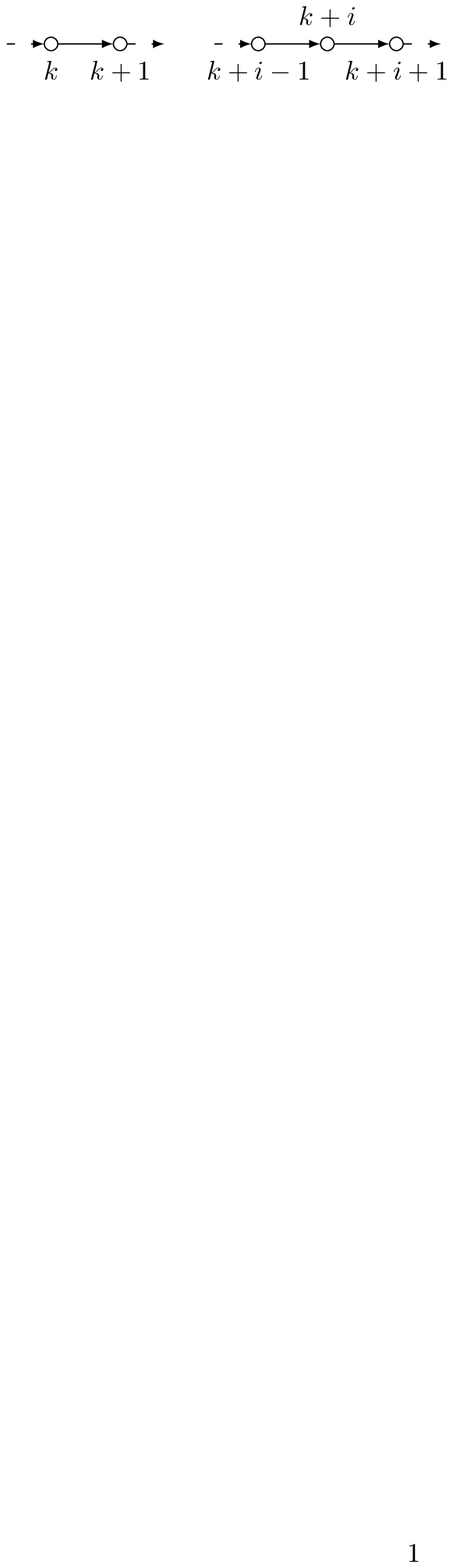}
\label{fig:linearNetwork}
\end{minipage}}
\subfloat[]{%
\begin{minipage}[c][0.5\width]{0.5\textwidth}
\centering%
\includegraphics{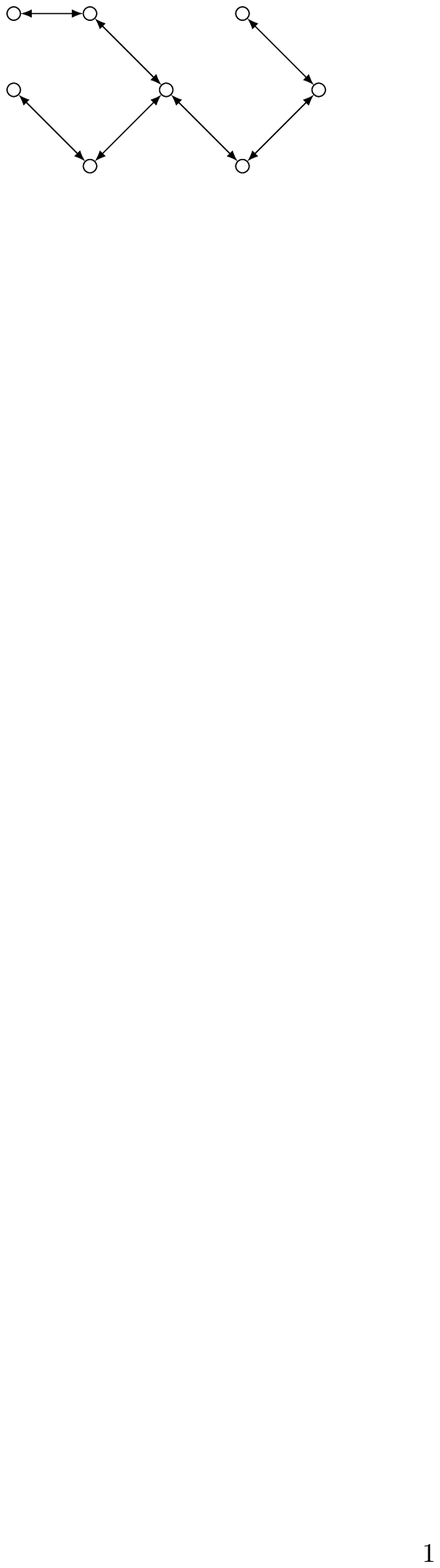}
\label{fig:treeNetwork}
\end{minipage}}
\caption{(a) A directed linear network. (b) A tree.}
\end{figure}

As a first structured example, consider a {\em directed linear network} with infinitely many nodes. See Figure~\ref{fig:linearNetwork}. In this situation we assume the channel is such that information from node $k$ can be directly transmitted only to node $k+1$. While channel interference may be taken into account, the model is insightful enough even in the case of perfect channel conditions. The choice that each node faces at any time instant is what information to transmit: its own or that of some node to the left of it. A Bernoulli policy is then determined by a probability distribution, $\{p_i,~ i=0,1,2,\ldots\}$ such that each node $k$ transmits or relays information about node $k-i$ with probability $p_i$.

For this class of networks, finding the marginal distribution of age is a
simple task. We assume stationarity and thus suppress the dependence on the time $n$. Denote by $A_{k+i,k}$ the age of information at some arbitrary node $k+i$ with respect to the information from node $k$. Then, for infinitely long networks, the random variables $A_{k+i,k}$ have the same distribution for every $k$, thus for shorthand we write $A_i$.   Now the time it takes information to propagate from node $k$ to node $k+i$ is distributed as the sum of $i$ independent geometric random variables (each with support $\{1,2,\ldots\}$) having parameters $p_0, p_1, \ldots,p_{i-1}$. Hence we have,
\[
\E[A_i] = \sum_{j=0}^{i-1} \frac{1}{p_j} \, ,
\qquad
\Var(A_i) = \sum_{j=0}^{i-1} \frac{1-p_j}{p_j^2} \, .
\]

A similar line of argumentation can be applied to infinite or finite trees as in Figure~\ref{fig:treeNetwork}. Since there is only one path\footnote{Throughout, we ignore redundant receptions in which a node receives information it has already relayed.} that information can take between any two nodes we again have that the set ${\{k: R_{k,i}(n) I_{k,j}(n) =1\}}$ appearing in \eqref{eqn:mainRecurs} contains at most one element. Thus the distribution of the age of information can be represented as a sum of independent geometric random variables (whose parameters depend generally both on the Bernoulli policy and on possible channel interference, in a straight-forward way). Further details are in~\cite{selen}.

\subsection*{Ring networks}
\begin{figure}
\centering
\includegraphics{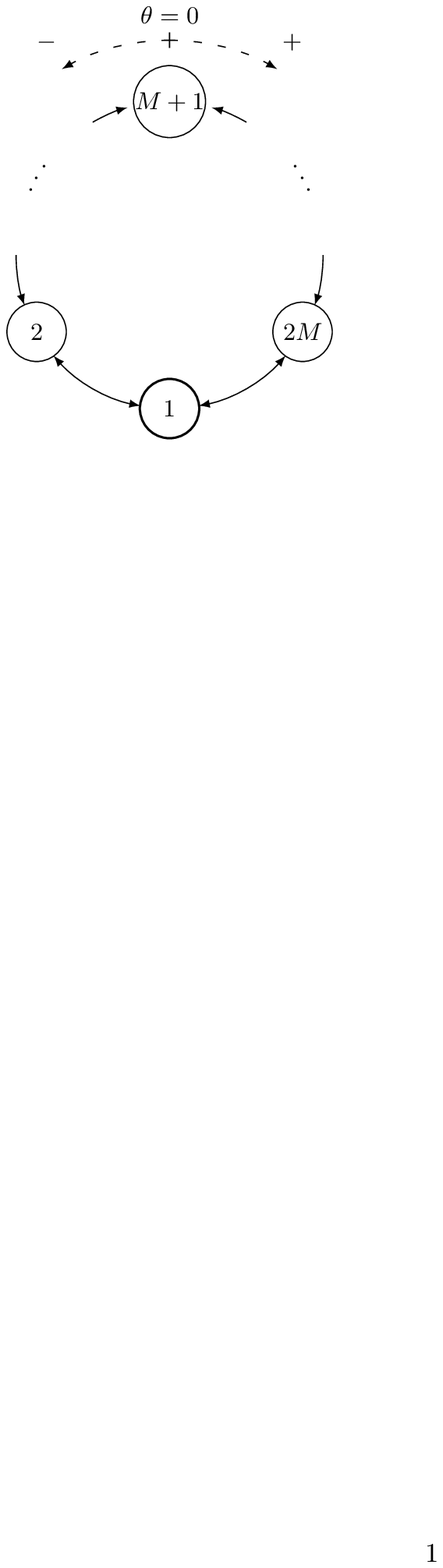}
\caption{Ring network of $2M$ nodes with node 1 the source of information.}
\label{fig:modelDescription}
\end{figure}

For modeling of situations in which information may travel on more than one route, a natural first step is to consider ring networks as in Figure~\ref{fig:modelDescription}. For brevity we consider networks with an even number of nodes, say $2M$ and assume ideal channels (a channel in which every transmitted packet is received). Each node transmits packets of information to its two closest neighbours. Assuming rotational symmetry, it is sufficient to study the distribution of the age of information with respect to a single source, say node 1. The age of information at node $i$ is then given by $A_{i,1}$, $i = 1,\ldots,2M$, for shorthand we write $A_i$.

Let us introduce a global coordinate variable $\theta$, defined for $i = 1,\ldots,2M$, by
\[
\theta:=\frac{i-1-M}{M} \in \{-M/M,(-M+1)/M,\ldots,0,\ldots,(M-1)/M\}.
\]
It is now convenient to use the value $Z_\theta := A_{M(\theta+1)+1}$. Figure~\ref{fig:modelDescription} illustrates both the node numbering $i$ and the coordinate variable $\theta$.

At every time slot, each node decides which sensor information it should relay (its own sensor information is also an option). Using Bernoulli policies, node $i$ transmits information it knows about node $j$ with a probability depending on the ``angle'' of $j$ relative to $i$, namely $\theta$. Denote this probability $q(\theta)$. 

Our aim is to study the marginal distribution $Z_\theta$. For each $\theta$, $Z_\theta$ is the minimum of the age of information coming from the clockwise and anticlockwise directions.  Information flowing back to the source is redundant, so these are equivalent to the age of information processes in a ring network with clockwise or anticlockwise {\em directed} transmission, respectively. Using the same reasoning as in the linear network, the age of information in one direction is distributed as a sum of independent geometric random variables in a ring network with directed transmission. We denote the directed age of information in the clockwise direction by $X_\theta^{(+)}$ and in the anticlockwise direction by $X_\theta^{(-)}$. To exemplify, $X_\theta^{(+)}$ is the age of information of the node corresponding to $\theta$ with respect to the source, node 1, when ignoring information coming from the anticlockwise direction.

Using the same reasoning as in the directed linear network, we note that
$X^{(+)}_{\theta}$ is a sum of independent geometric random variables, with
\[
\E[X^{(+)}_{\theta}] = \sum_{d=-M}^{\theta M-1} \frac{1}{q(d/M)} \, ,
\qquad
\Var(X^{(+)}_{\theta}) = \sum_{d=-M}^{\theta M-1} \frac{1-q(d/M)}{q(d/M)^2} \, .
\]
In the anticlockwise directed transmission ($X_\theta^{(-)}$) and with $q(\cdot)$ symmetric with respect to the distance from the source, the mean and variance are expressed in the same way except for the interchange of $\theta$ by $(1-\theta)$ in the summation.

As a Bernoulli policy, we suggest a parametric family of distributions:
\begin{align*}
q(\theta) &:= \left\{ \begin{array}{ll}
\beta, & \theta = -1 \, , \\
C\alpha^{|M(|\theta|-1)|}, & \theta > -1 \, , \\
\end{array} \right. ~~\mbox{where}~~
C := \left\{ \begin{array}{ll}
\frac{1-\beta}{2M-1}, & \alpha = 1 \, , \\
\frac{(1-\beta)(1-\alpha)}{2\alpha-\alpha^M(\alpha+1)}, & \alpha < 1 \, , \\
\end{array} \right.
\end{align*}
where $\alpha \in (0,1]$ describes the geometric decay in probability when moving away from the source and $\beta \in (0,1)$ is the probability mass of the source transmitting its own information. 

This family allows various behaviours: A uniform transmission probability ($\alpha=1$) or alternatively decaying probabilities when moving further away from the source ($\alpha<1$), both with or without a different probability of transmitting at the source as determined by $\beta$.  The information sent by the source is usefully transmitted in both the clockwise and anticlockwise direction, whereas relayed information only benefits one of the relay's neighbours. This suggests that $\beta$ should give a higher weight to the source; we optimize $\beta$ in Section~\ref{sec:asymptoticsAndOptimPolicy}.

\section{Explicit and algorithmic solutions}
\label{sec:explicit}
Finding the stationary distributions, their marginals or the means of our models is in general not straightforward. Nevertheless in this section we report some successful results. In doing so we illustrate a recurring pattern in these types of models: Using marginal distributions to find joint distributions.

The most basic model is a sensor node transmitting to a receiver, where there is a chance of $\lambda \in(0,1)$ for successful reception.  In this case the age of information at the receiver follows a specific GI/M/1 type Markov chain (c.f. \cite{asmussen}, Section XI.3) in which transitions increment the state by one with probability $(1-\lambda)$ or reset the state to 0 with probability $\lambda$.  As with all GI/M/1 (scalar) Markov chains, the stationary distribution is geometric, in this case with parameter $(1-\lambda)$ and support $\{0,1,\ldots\}$. We shift the support to $\{1,2,\ldots\}$ to accommodate the minimal possible age, $1$. In general the value of $\lambda$ may be influenced by both the channel properties and the transmission policy. For example we may have $\lambda = p q$ where $p$ is the chance of receiving a packet conditional on it being transmitted and $q$ is the chance of transmitting.

This GI/M/1 type stationary distribution can be used as a ``building block'' for finding the (multi-dimensional) stationary distributions of more complicated models. We illustrate this now for two types of models: star networks and a small ring, further examples and details are in \cite{selen}.

\subsection*{Star networks}
\begin{figure}
\centering
\includegraphics{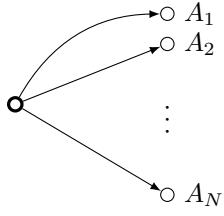}
\caption{A star topology.
}
\label{fig:singleHops}
\end{figure}

Consider star networks as illustrated in Figure~\ref{fig:singleHops}. Transmissions take place from the source node to $N$ receivers. We denote a version of the steady state age of information at node $i$ with respect to the source node by $A_i$.  What is then the joint distribution of $A_1,\ldots,A_N$?

\begin{figure}
\centering
\includegraphics{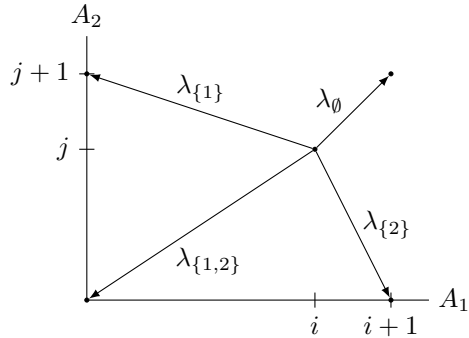}
\caption{Markov chain transition diagram for a star with $N=2$.}
\label{fig:stateSpace2Hops}
\end{figure}

To illustrate the solution approach we first consider the case of $N=2$. Let $\lambda_{\emptyset}$, $\lambda_{\{1\}}$, $\lambda_{\{2\}}$ and $\lambda_{\{1,2\}}$ denote the respective probabilities that reception occurs at neither node, node $A_1$ only, node $A_2$ only, or both nodes. The transition diagram of this model is shown in Figure~\ref{fig:stateSpace2Hops}.

Let $\pi_{i,j} :=\Prob(A_1 = i,A_2 = j)$. Then,
\begin{subequations}
\begin{align}
\pi_{1,1} &= \lambda_{\{1,2\}} \sum_{i=1}^\infty \sum_{j=1}^\infty \pi_{i,j} \, , \label{eqn:2TR_1RE_pi0} \\
\pi_{i,1} &= \lambda_{\{2\}} \sum_{j=1}^\infty \pi_{(i-1),j}, \quad i \geq 2 \, , \label{eqn:2TR_1RE_2} \\
\pi_{1,j} &= \lambda_{\{1\}} \sum_{i=1}^\infty \pi_{i,(j-1)}, \quad j \geq 2 \, , \label{eqn:2TR_1RE_3} \\
\pi_{i,j} &= \lambda_{\emptyset} \cdot \pi_{(i-1),(j-1)}, \quad i,j \geq 2 \, .
\end{align}
\end{subequations}
%

Now a key observation is that in \eqref{eqn:2TR_1RE_2}-\eqref{eqn:2TR_1RE_3}
there is summation over one entire coordinate, therefore we can use the
marginal distributions. For nodes $k=1,2$, let $c_k =
1-(\lambda_{\{3-k\}}+\lambda_{\{1,2\}})$ denote the probability of no
reception on the other node, $3-k$. Then as in the GI/M/1 type Markov chain described above, the marginal distributions are given by
\[
\pi_i^{(A_k)} :=
 \Prob(A_k = i) = (1-c_{3-k})c_{3-k}^{i-1}, \quad k = 1,2, ~ i = 1,2,\ldots  \> .
\]
Since $\pi_{1,1} = \lambda_{\{1,2\}}$, these marginal distributions imply the equilibrium equations simplify to
\begin{align*}
\pi_{1,1} &= \lambda_{\{1,2\}} \, ,\\
\pi_{i,1} &= \lambda_{\{2\}} \pi_{i-1}^{(A_1)} = \lambda_{\{2\}} c_2^{i-2} (1-c_2), \quad i \geq 2 \, ,   \\
\pi_{1,j} &= \lambda_{\{1\}} \pi_{j-1}^{(A_2)} = \lambda_{\{1\}} c_1^{j-2} (1-c_1), \quad j \geq 2 \, ,  \\
\pi_{i,j} &= \lambda_{\emptyset} \cdot \pi_{(i-1),(j-1)}, \quad i,j \geq 2 \,.
\end{align*}
These then yield the stationary distribution
\begin{align*}
\pi_{i,j} = \left\{ \begin{array}{cc} \lambda_{\emptyset}^{i-1} \lambda_{\{1,2\}}, & i=j \, , \\
                                      \lambda_{\emptyset}^{j-1} \lambda_{\{2\}} c_2^{i-j-1} (1-c_2), & i>j \, , \\
                                      \lambda_{\emptyset}^{i-1} \lambda_{\{1\}} c_1^{j-i-1} (1-c_1), & i<j \, .
            \end{array} \right.
\end{align*}
After some straightforward calculations this yields
\begin{align*}
\mathrm{Cov}(A_1,A_2) &= \frac{ \lambda_{\emptyset}\lambda_{\{1,2\}} - \lambda_{\{1\}}\lambda_{\{2\}} }{(\lambda_{\{1\}} + \lambda_{\{1,2\}})(\lambda_{\{2\}} + \lambda_{\{1,2\}})(1-\lambda_{\emptyset})} \, .
\end{align*}
It can now be verified that if there is no interaction between the
communication links, i.e.,
$(\lambda_{\{1\}}+\lambda_{\{1,2\}})(\lambda_{\{2\}}+\lambda_{\{1,2\}}) =
\lambda_{\{1,2\}}$, then there is a product form solution to $\pi_{i,j}$ and the covariance is $0$. Otherwise, the covariance is non-zero and can be used to get LMMSE (linear minimum mean squared error estimates) of $A_k$ based on $A_{3-k}$. We do not discuss this further here.

The idea of a network with $N=2$ can now be generalized to arbitrary $N$ by recursive usage of marginal distributions of some lower order. We describe this in brief and present an algorithm for calculating the {\em exact} stationary distribution.

The Bernoulli policies and i.i.d. channel conditions imply that we may essentially have $\lambda_B$ for any receiving subset of nodes $B$. We let $D$ denote some proper subset of the set of all nodes in order to consider smaller networks in the recursive specification that follows. When computing the joint distribution of a subset $D$ of all the nodes, we need to know $\lambda_{B;D}$, which are the probabilities of successful reception on the nodes in the set $B$, given that we only consider receptions on the nodes in the subset $D$ of the full network. That is, we ignore transmissions to the nodes in the complement of $D$.

Let $D = \{i_1,i_2,\ldots,i_{|D|}\}$. To find
\[
\pi_{a_{i_1}, \ldots,a_{i_{|D|}}}^{(D)} := \Prob(A_{i_1} = a_{i_1}, \ldots, A_{i_{|D|}} = a_{i_{|D|}}) \, ,
\]
let $j_1,j_2,\ldots,j_{|D|}$ be a permutation of $D$ such that $0 \le a_{j_1} \le a_{j_2} \le \ldots \le a_{j_{|D|}}$, and let $\tilde \pi^{(D)}_{a_{j_1},a_{j_2},\ldots,a_{j_{|D|}}} = \pi^{(D)}_{a_{i_1}, a_{i_2},\ldots,a_{i_{|D|}}}$.
Then $\tilde \pi$ can be calculated recursively by Algorithm~\ref{algo:jointSingleHops}.

\begin{algorithm}
\caption{Joint distribution of $|D|$ nodes in a network of $N$ nodes.}
\label{algo:jointSingleHops}
\begin{algorithmic}[1]
\If{$a_{j_1} =  \ldots = a_{j_{|D|}}$}
    \State $m=|D|$
\Else
    \State $m = \min \{ k : a_{j_k} < a_{j_{k+1}}\}$
\EndIf
\If{$a_{j_{|D|}} = 1$}
    \State $\tilde \pi^{(D)}_{1,\ldots,a_{j_{|D|}}} = \lambda_{D;D}$ \label{eqn:singleHops_origin}
\Else
    \State Let $F = D\setminus\{j_{1},\ldots,j_{m}\}$ and its complement is $F^c$. 
    \State
    \vspace{-15pt}
        \begin{align*}
        \tilde \pi^{(D)}_{a_{j_1},a_{j_2},\ldots,a_{j_{|D|}}} &= \left\{ \begin{array}{ll}
        \lambda_{\emptyset;D}^{a_{j_1}-1} \cdot \tilde \pi^{(D)}_{1,\ldots,1,a_{j_{m+1}}-a_{j_1}+1,\ldots,a_{j_{{|D|}}}-a_{j_1}+1}, & a_{j_1} > 1 \\
        \lambda_{F^c;C} \cdot \tilde \pi^{(F)}_{a_{j_{m+1}}-1,\ldots,a_{j_{|D|}}-1}, &  a_{j_1} = 1
        \end{array} \right.
        \end{align*} \label{eqn:singleHops_rest}
    \vspace{-15pt}
\EndIf
\end{algorithmic}
\end{algorithm}

Similarly to the $N=2$ case, the probability at the point $(1,1,\ldots,1)$
equals the probability of reception on all nodes in the network; see line
\ref{eqn:singleHops_origin}. For any state in the interior of the state space,
i.e., the smallest age satisfies $a_{i_1} > 1$, we can compute the probability by moving back along the diagonal to the nearest (hyper) plane or edge and using the knowledge that there is a geometric decay along the diagonals. This is shown in the first part of the equation on line \ref{eqn:singleHops_rest}. If we are already on a (hyper) plane or edge, we can use the marginal distribution of all the other nodes that have a strictly positive age (see second part of line \ref{eqn:singleHops_rest}). See report~\cite{selen} for an illustration in the case of $N=3$.

\pagebreak
\subsection*{A small ring}
\begin{figure}[H]
\centering
\includegraphics{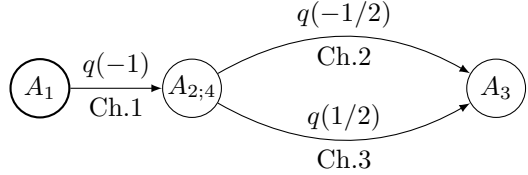}
\caption{Alternate representation of a network with four nodes where node 1 is the source. There are three channels.}
\label{fig:explicitM4}
\end{figure}

Let us now consider the smallest non-trivial ring: a ring with $2M = 4$ nodes. We exploit now the fact that $A_2 = A_4$ and denote them by $A_{2;4}$. We then allow this ``virtual'' node $A_{2;4}$ to transmit over two separate channels to the node diametrically opposite the source. Thus we can represent the steady state age of information by $A_3 \in \{2,3,\ldots\}$ and $A_{2;4} \in \{1,2,\ldots,A_3\}$. See Figure~\ref{fig:explicitM4}.

Denote $\pi_{i,j} := \Prob(A_{2;4} = i, A_3 = j)$. Observe that the marginal
distribution of $A_{2;4}$ is geometric with parameter
$q(-1)=\beta$ and support $\{1,2,\ldots\}$, as we found earlier in this
section. Let $\pi_i^{(A_{2;4})}=\Prob(A_{2;4}=i)$. Similarly to the star, this
value appears in the balance equations of $\pi_{i,j}$. These equations are
based on reception probabilities on subsets of the channels denoted by
$\lambda_B$, where $B$ is a set of channels. For example $\lambda_{\{2,3\}} =
(1-q(-1))q(-1/2)q(1/2)$.
\begin{subequations}
\begin{align}
\pi_{1,2} &= \left( \lambda_{\{1,2\}} + \lambda_{\{1,3\}} + \lambda_{\{1,2,3\}} \right) \pi_1^{(A_{2;4})} \, , \label{eqn:beginExactEq}\\
\pi_{2,2} &= \left( \lambda_{\{2\}} + \lambda_{\{3\}} + \lambda_{\{2,3\}}\right) \pi_1^{(A_{2;4})} \, , \\
\pi_{1,j} &=  \displaystyle \lambda_{\{1\}} \sum_{i=1}^{j-1} \pi_{i,j-1}
 +\left( \lambda_{\{1,2\}} + \lambda_{\{1,3\}} + \lambda_{\{1,2,3\}} \right) \pi_{j-1}^{(A_{2;4})}, \quad j \geq 3  \, , \label{eqn:beginAlg} \\
\pi_{i,j} &= \left( \lambda_{\{2\}} + \lambda_{\{3\}} + \lambda_{\{2,3\}}\right)\pi_{i-1}^{(A_{2;4})}  + \lambda_{\emptyset}\pi_{i-1,j-1}, \quad i = j, ~ i \geq 3 \, , \label{eqn:midAlg}\\
\pi_{i,j} &= \lambda_{\emptyset}\pi_{i-1,j-1}, \quad i \neq j, ~ i \geq 2, ~ j \geq 3 \, . \label{eqn:endAlg}
\end{align}
\end{subequations}
Algorithm~\ref{algo:jointM4} uses these equations to calculate $\{\pi_{i,j},~
i,j \le K\}$ {\em exactly} for any $K$.
\begin{algorithm}[h]
\caption{Joint distribution of $A_{2;4}$ and $A_3$.}
\label{algo:jointM4}
\begin{algorithmic}
\State Use the known $\pi_{i,j}^{(A_{2;4})}$ and set $\pi_{1,2}$ and $\pi_{2,2}$.
\For{$j=3:K$}\Comment{Iterate until a bounding box of size $K$ is reached.}
    \For{$i=1:j$}
        \If{$i = 1$}
            \State Calculate $\pi_{i,j}$ using \eqref{eqn:beginAlg}, based on $\pi_{j-1}^{(A_{2;4})}$ and $\pi_{i,j-1}$.
        \ElsIf{$i=j$}
            \State Calculate $\pi_{i,j}$ using \eqref{eqn:midAlg}, based on $\pi_{i-1}^{(A_{2;4})}$ and $\pi_{i-1,j-1}$.
        \Else
            \State Calculate $\pi_{i,j}$ using \eqref{eqn:endAlg}, based on $\pi_{i-1,j-1}$.
        \EndIf
    \EndFor
\EndFor
\end{algorithmic}
\end{algorithm}

\begin{figure}[h]
\centering
\subfloat[$(\alpha=1,\beta=\frac{1}{2M})$ Bernoulli policy]{
\includegraphics{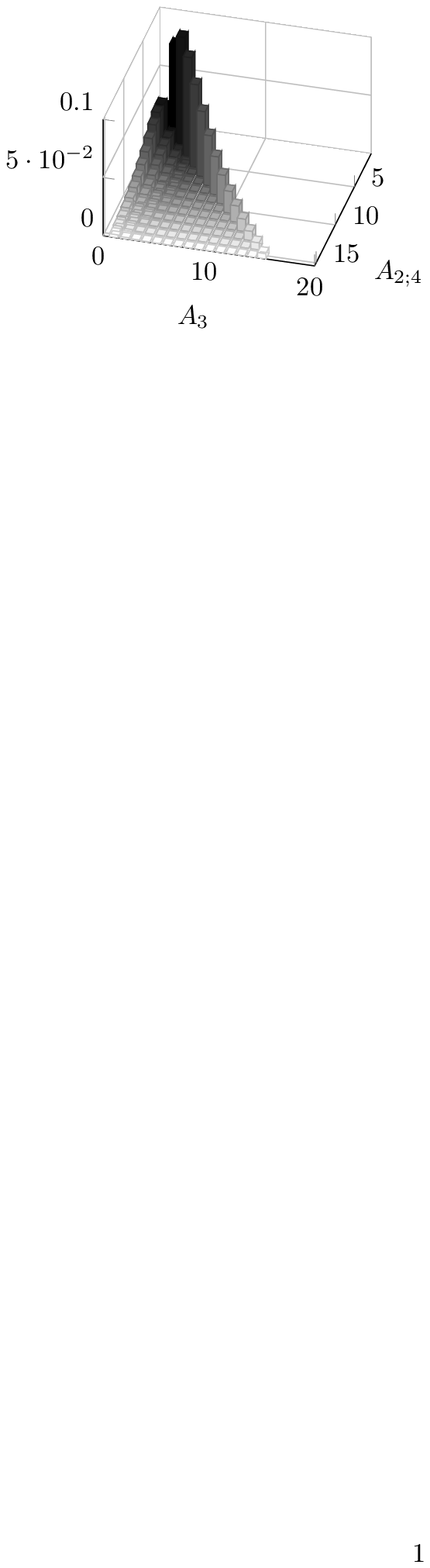}
}
\hspace{2mm}
\subfloat[$(\alpha=0.1,\beta=\frac{2}{2M})$ Bernoulli policy]{
\includegraphics{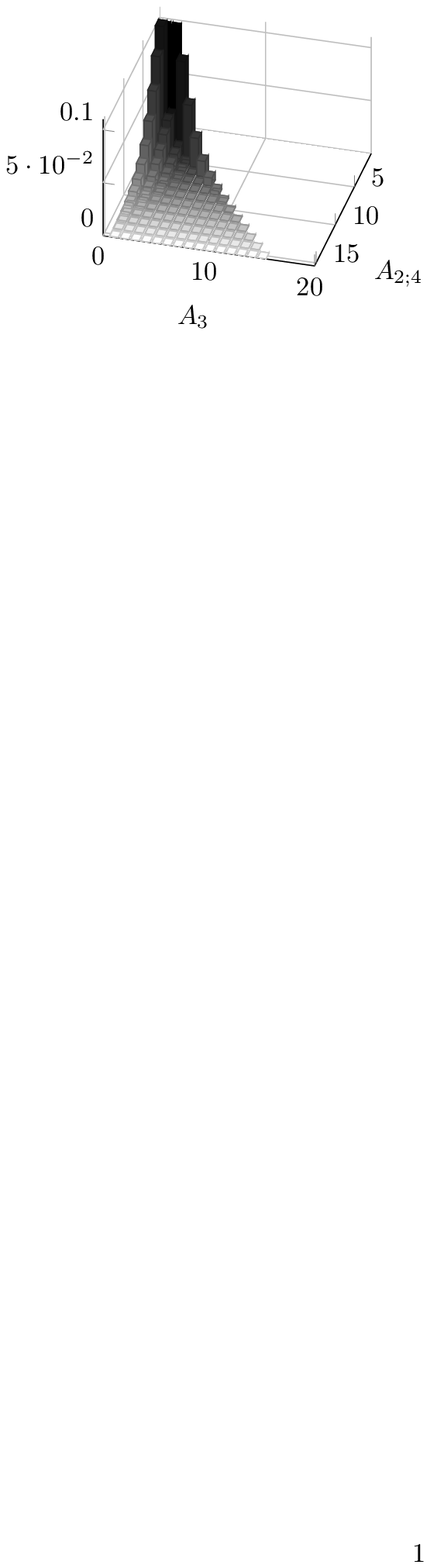}
}
\caption{Joint distribution of $A_{2;4}$ and $A_3$ for two different policies, using Algorithm~\ref{algo:jointM4}.}
\label{fig:surfaceSmallNetwork}
\end{figure}

We now present a numerical example. We compute the joint distribution of $A_{2;4}$ and $A_3$ for two sets of transmission parameters $(\alpha,\beta)$. The first is a uniform policy $(\alpha=1,\beta=\frac{1}{2M})$, and the second has its probability mass concentrated around the source, $(\alpha=0.1,\beta=\frac{2}{2M})$. Figure~\ref{fig:surfaceSmallNetwork} shows the joint distribution found by Algorithm~\ref{algo:jointM4}. In the first case the probability mass of the joint distribution is more widely spread out over the state space and in the latter it is more concentrated around the minimum ages, i.e.~$a_{2;4} = 1$ and $a_3 = 2$.

A similar approach to that of Algorithm~\ref{algo:jointM4} can essentially be applied to networks with more nodes. However, this is analytically demanding and becomes impractical. Even for a network with 5 nodes there are 5 possible transmissions and thus $2^5$ different subsets of $B$ in $\lambda_B$ and many more equations in comparison to \eqref{eqn:beginExactEq}-\eqref{eqn:endAlg}. We therefore shift our attention to approximations.

\section{Asymptotic approximations in rings}
\label{sec:asymptoticsAndOptimPolicy}
In this section we present an asymptotic evaluation of ring networks with
$\alpha=1$ and some $\beta$. We revisit the question presented in the
introduction: How should a node balance transmitting its own information against relaying? Alternatively, what is a good value for $\beta$? Our analysis is based on the representation
\[
Z_{\theta} = \tilde{X}_{\frac{M-1}{M}} +
\big(\tilde{X}^{(+)}_{\theta}
\wedge
\tilde{X}^{(-)}_{\theta}\big) \, ,
\]
where $\tilde{X}_{\frac{M-1}{M}}$ represents the age at the neighbouring nodes
of the source (both have the same age) and $\tilde{X}_{\theta}^{(+)},
\tilde{X}_{\theta}^{(-)}$ represent the age difference between the node in question and the neighbouring nodes of the source, in the clockwise and anticlockwise directions respectively, based on directed transmission.

For large $M$, we are guided by the central limit theorem to use a Gaussian
approximation for each of the directed transmissions, i.e., the Negative Binomially distributed $\tilde{X}^{(+)}_{\theta}$ and $\tilde{X}^{(-)}_{\theta}$ are approximately normally distributed with
\[
\mu_\theta^{(+)}:=\E[\tilde{X}^{(+)}_{\theta}] = ((1+\theta)M-1) C^{-1} \, ,
\qquad
\mu_\theta^{(-)}:=\E[\tilde{X}^{(-)}_{\theta}] = ((1-\theta)M-1) C^{-1} \, ,
\]
and standard deviations, $\sigma^{(+)}_{\theta} :=\sqrt{\mu^{(+)}_{\theta}\left(C^{-1}-1\right)}$,
$\sigma^{(-)}_{\theta}:=\sqrt{\mu^{(-)}_{\theta}\left(C^{-1}-1\right)}$ respectively.  We now have
\[
Z_{\theta} \approx^d  \hat{Z}_\theta := \tilde{X}_{\frac{M-1}{M}} + (\mathcal{N}_{\tilde{X}^{(+)}_{\theta}} \wedge\mathcal{N}_{\tilde{X}^{(-)}_{\theta}}) \, ,
\]
where $\approx^d$ informally denotes approximate equality in distribution and the $\mathcal{N}$ variables are independent versions of normal random variables with the aforementioned parameters. In this paper we do not formalize this as a weak-convergence result (as $M \to \infty$).  This technical hurdle is left for future research.

In \cite{hunter} (see also \cite{tong}) the moments of the minima of normally distributed random variables are given.  We exploit these results here to find approximating expressions for the mean and variance of $Z_\theta$. Denoting the CDF and PDF of the standard normal distribution by $\Phi(\cdot)$ and $\phi(\cdot)$ respectively,  we obtain
\begin{align}
\E[\hat{Z}_{\theta}] &= \frac{1}{\beta}+ \mu^{(+)}_{\theta} \Phi\left(\frac{-\bar{\mu}}{\Delta}\right) + \mu^{(-)}_{\theta} \Phi\left(\frac{\bar{\mu}}{\Delta}\right) - \Delta \phi\left(\frac{-\bar{\mu}}{\Delta}\right) \, , \label{eqn:approxMean}\\
\E[\hat{Z}_{\theta}^2] &= \frac{-1}{\beta}+ \omega^{(+)}\Phi\left(\frac{-\bar{\mu}}{\Delta}\right) + \omega^{(-)}\Phi\left(\frac{\bar{\mu}}{\Delta}\right) - \bar{\mu} \Delta \phi\left(\frac{-\bar{\mu}}{\Delta}\right) \, , \label{eqn:approxVar}
\end{align}
where $\bar{\mu} := \mu^{(+)}_{\theta} -\mu^{(-)}_{\theta}$, $\Delta := \sqrt{\left(\sigma^{(+)}_{\theta}\right)^2+\left(\sigma^{(-)}_{\theta}\right)^2}$, $\omega^{(+)} := \left(\mu^{(+)}_{\theta}\right)^2+\left(\sigma^{(+)}_{\theta}\right)^2$  and  $\omega^{(-)}:=\left(\mu^{(-)}_{\theta}\right)^2+\left(\sigma^{(-)}_{\theta}\right)^2$.

\begin{figure}
\centering
\subfloat{
\includegraphics{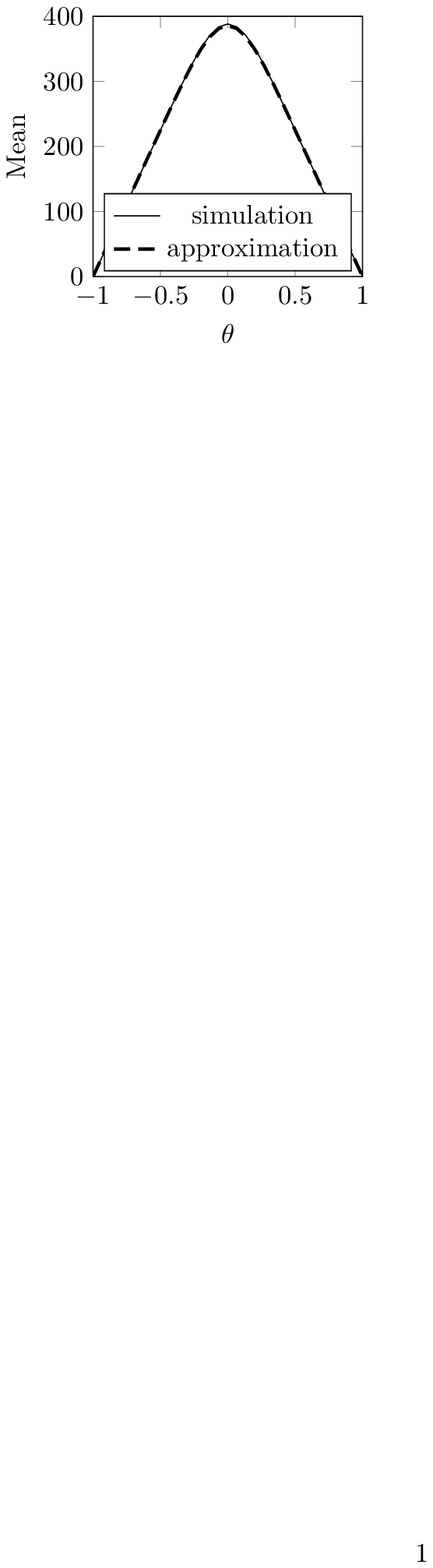}
}
\hspace{2mm}
\subfloat{
\includegraphics{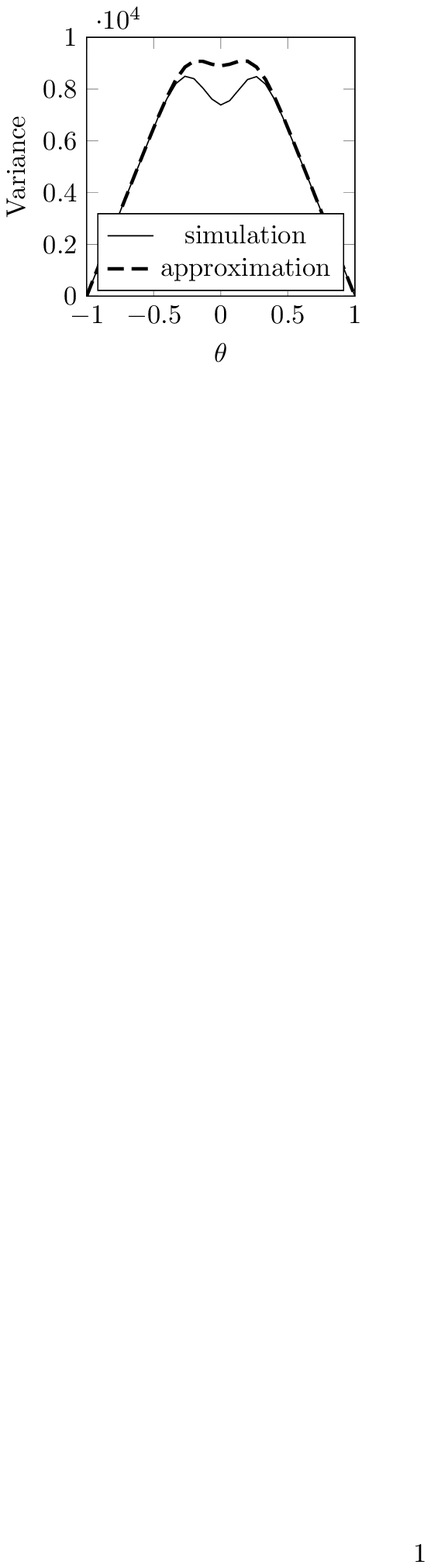}
}
\caption{Comparison of the mean and variance obtained from simulation with the approximation values $\E[\hat{Z}_\theta]$ and $\Var[\hat{Z}_\theta]$. This is for $2M=30$.
}
\label{fig:simulationApproximation30}
\end{figure}

We conjecture that for any $\theta \in [-1,1]$, $\lim_{M \to \infty} \E[Z_{[\theta]}]/\E[\hat{Z}_\theta]=1$ and the same for the variance (here $[\theta]$ denotes the nearest value that $\theta$ may attain over the grid). We have verified this conjecture numerically by means of extensive Monte-Carlo simulations.  As an illustration we compare the curves for $2M=30$ nodes in Figure~\ref{fig:simulationApproximation30}.

To observe convergence of the variance to the ``volcano curve'', we also ran long simulations for $2M=100$. The results are not displayed here. Our numerical experiments have also clearly indicated that $\E[Z_\theta] \geq \E[\hat{Z}_\theta]$ and $\Var(Z_\theta) \leq \Var(\hat{Z_\theta})$. We leave proofs of these inequalities for future work.

\subsection*{The asymptotically best $\beta$}
\begin{figure}[h]
\centering
\includegraphics{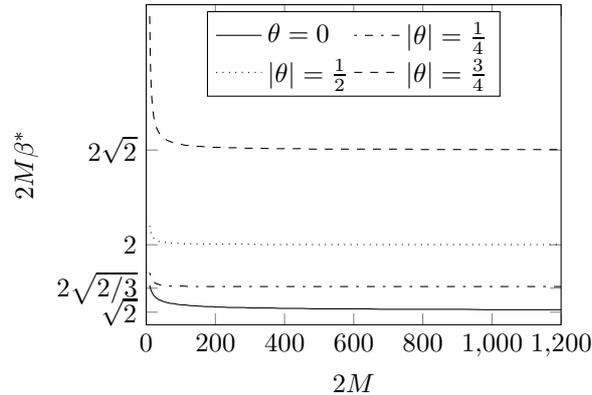}
\caption{Optimal values for $\beta$ for various angles $\theta$ for increasing network size $2M$.
}
\label{fig:optimalBeta}
\end{figure}

We now optimize the transmission policy with respect to minimizing the mean age of information at the node corresponding to $\theta$. For $\theta = 0$ we can simplify \eqref{eqn:approxMean}. We know that the mean and variance of $\tilde{X}^{(+)}_0$ and $\tilde{X}^{(- )}_0$ are equal and we omit the superscripts $(+), (-)$. This leads to the following expression:
\[
\E[\hat{Z}_0] = \frac{1}{\beta} + \mu_0 - \frac{\sigma_0}{\sqrt{\pi}} \, .
\]
The mean $\mu_0$ is $O(M^2)$, whereas $\sigma_0$ is $O(\sqrt{M^3})$ and both
scale with $1/(1-\beta)$. Hence the mean dominates the standard deviation for
large $M$, and thus
\[
\E[\hat{Z}_0] \approx
\frac{1}{\beta} + \frac{2M^2}{1-\beta}
\]
for large $M$. This is minimized for
$\hat{\beta}^* = \sqrt{2}/(2M)$ for large $M$ and $\theta = 0$. For $\theta\neq 0$,
again $\sigma_{\theta} = o(\mu_{\theta})$ whence
$|\bar\mu/\Delta|\rightarrow\infty$ and $\hat{\beta}^* =
\sqrt{\frac{2}{1-|\theta|}}/(2M)$. We numerically compute the $\beta^*$ values for various fractions $\theta$ using \eqref{eqn:approxMean}, summarized in Figure~\ref{fig:optimalBeta}. Observe the converge of $\beta^*$ to $\hat{\beta}^*$ as $M \to \infty$.

We have thus found that in large rings, if the overall goal is to maintain timely information at the farthest node from each sensor, then each node should transmit its own information about $40\%$ more frequently than the information of some other node. This finding is of course based on a series of assumptions and stylized modeling assumptions.  Yet it can perhaps serve as a rule of thumb for gossip networks, if there is no better alternative.

\section{Conclusion} \label{sec:conclusion}
We have developed a simple Markovian framework for the design and analysis of a
gossip protocol in tree or ring topology networks where information is
probabilistically updated by each individual node and sent over
bandwidth-limited lossy wireless channels. Using the framework, we presented
some basic results for the mean, the variance, and distribution of age in the
studied star networks and a small ring, including non-trivial explicit and algorithmic solutions to obtain the age of information distribution. For large ring networks, we obtained asymptotic forms for the age of information using normal approximations and explored the optimal strategy to forward information in such a network.

Future work will deal with the extension of the framework beyond the linear (or
tree) and ring network topologies where new asymptotic approximations could be
developed. For most applications, including ITS, information about nearby nodes
is more important than information about distant nodes.  Hence it will be
useful to consider both optimizing weighted means of ages of information, and
also coarsely aggregating information as it emanates further from its source.
We also wish to settle some of the conjectures laid out regarding the ring asymptotics.

\subsection*{Acknowledgment} \label{sec:ack}
This research began while the first author was visiting Swinburne University of Technology and The University of Queensland and also took place while the first and second author were visiting the University of Haifa.
This work was supported by Australian Research Council (ARC) grants DP130100156, DE130100291, FT0991594 and  FT120100723. The authors thank Ivo Adan for useful comments.


\end{document}